\documentclass[11pt]{amsart}

\usepackage{amsmath, amsfonts, amssymb}
\usepackage{graphics}
\usepackage[usenames]{color}

{\theoremstyle{plain}
\newtheorem{theorem}{Theorem}[section]

\newtheorem{proposition}[theorem]{Proposition}
\newtheorem{lemma}[theorem]{Lemma}
\newtheorem{corollary}[theorem]{Corollary}
}

{\theoremstyle{definition}

\newtheorem{remark}[theorem]{Remark}

}

\newcommand{\Sn}{\mathfrak{S}_n}
\newcommand{\Dn}{\mathfrak{D}_n}
\newcommand{\Dnr}{\mathfrak{D}_n^{(r)}}
\newcommand{\Cr}{\mathcal{C}_r}
\newcommand{\CS}{\mathcal{C}_r \wr \mathfrak{S}_n}

\newcommand{\Des}{\mathrm{Des}}
\newcommand{\maj}{\mathrm{maj}}
\newcommand{\des}{\mathrm{des}}
\newcommand{\exc}{\mathrm{exc}}
\newcommand{\sub}{\mathrm{sub}}
\newcommand{\sign}{\mathrm{sgn}}

\newcommand{\der}{\mathrm{dp}}

\begin{document}
\title[Descent cynical germ]{Cyclic derangements}

\author{Sami H. Assaf}
\address{Department of Mathematics, Massachusetts Institute of Technology, 77 Massachusetts Avenue, Cambridge, MA 02139, USA}
\email{sassaf@math.mit.edu}
\thanks{Partially supported by NSF Mathematical Sciences Postdoctoral
  Research Fellowship DMS-0703567}



\begin{abstract} 
  A classic problem in enumerative combinatorics is to count the
  number of derangements, that is, permutations with no fixed
  point. Inspired by a recent generalization to facet derangements of
  the hypercube by Gordon and McMahon, we generalize this problem to
  enumerating derangements in the wreath product of any finite cyclic
  group with the symmetric group. We also give $q$- and
  $(q,t)$-analogs for cyclic derangements, generalizing results of
  Brenti and Gessel.
\end{abstract}

\maketitle

\section{Derangements}\label{sec:derangements} 

A derangement of $\{1,2,\ldots,n\}$ is a permutation that leaves no
letter fixed. Algebraically, this is an element $\sigma$ of the
symmetric group $\Sn$ such that $\sigma(i) \neq i$ for any $i$, or,
equivalently, no cycle of $\sigma$ has length $1$. Geometrically, a
derangement of $\{1,2,\ldots,n\}$ is an isometry in $\mathbb{R}^{n-1}$
of the regular $(n-1)$-simplex that leaves no facet
unmoved. Combinatorially, these are matrices with entries from
$\{0,1\}$ such that each row and each column has exactly one nonzero
entry and no diagonal entry is equal to $1$. 

Let $\Dn$ denote the set of derangements in $\Sn$, and let $d_n =
|\Dn|$. The problem of enumerating derangements is the canonical
example of the principle of Inclusion-Exclusion \cite{EC1}:
\begin{equation}
  d_n = n! \sum_{i=0}^{n} \frac{(-1)^i}{i!}.
\label{e:derange}
\end{equation}
For example, the first few derangement numbers are $1, 0, 1, 2, 9, 44,
265$.

From (\ref{e:derange}) one can immediately compute that the
probability that a random permutation has no fixed points is
approximately (and very nearly) $1/e$. Another exercise that often
accompanies counting derangements is to prove the following recurrence
relation for $n \geq 2$,
\begin{equation}
  d_n = (n-1) \left( d_{n-1} + d_{n-2} \right),
\label{e:recur}
\end{equation}
with initial conditions $d_0=1$ and $d_1=0$; see \cite{EC1}. From
\eqref{e:recur} one can derive the following single term recurrence
for derangement numbers,
\begin{equation}
  d_n = n d_{n-1} + (-1)^n.
\label{e:recur1}
\end{equation}

Recently, Gordon and McMahon \cite{GoMc2009} looked at the problem of
enumerating isometries of the $n$-dimensional hypercube that leave no
facet unmoved. Algebraically, such an isometry is an element $\sigma$
of the hyperoctahedral group $B_n$ for which $\sigma(i) \neq i$ for
any $i$. Combinatorially, the problem then is to enumerate $n \times
n$ matrices with entries from $\{0,\pm 1\}$ such that each row and
column has exactly one nonzero entry and no diagonal entry equals
$1$. Using the same technique of Inclusion-Exclusion, Gordon and
McMahon derive a formula for the number of facet derangements similar
to \eqref{e:derange}, an expression of facet derangements in terms of
permutation derangements, and recurrence relations for facet
derangements similar to \eqref{e:recur} and \eqref{e:recur1}.

In Section~\ref{s:cyclic}, we consider elements $\sigma$ in the wreath
product $\CS$, where $\Cr$ is the finite cyclic group of order $r$ and
$\Sn$ is the symmetric group on $n$ objects. A \textit{cyclic
  derangement} is an element of $\CS$ with no fixed point. Denote the
set of cyclic derangements of $\CS$ by $\Dnr$, and denote their number
by $d_n^{(r)} = |\Dnr|$. Combinatorially, $d_n^{(r)}$ is also the
number of matrices with entries from
$\{0,1,\zeta,\ldots,\zeta^{r-1}\}$, where $\zeta$ is a primitive $r$th
root of unity, such that each row and each column have exactly one
nonzero entry with no diagonal entry equal to $1$. Using
Inclusion-Exclusion, we derive a formula for $d_n^{(r)}$ that
specializes to \eqref{e:derange} when $r=1$ and to the Gordon-McMahon
formula for facet derangements when $r=2$. We also give an expression
for $d_n^{(r)}$ in terms of permutation derangements as well as a two
recurrence relations specializing to \eqref{e:recur} and
\eqref{e:recur1} when $r=1$.

Another direction for generalizations is to $q$-count derangements by
various statistics. Gessel \cite{GeRe1993} introduced a $q$-analog for
derangements of $\Sn$, $q$-counted by the major index, that has
applications to character theory \cite{Schocker2003}. In
Section~\ref{s:qtcyclic}, we give a $q,t$-analog for cyclic
derangements of $\CS$ $q$-counted by a generalization of major index
and $t$-counted by signs that specializes to Gessel's formula at $r=1$
and $t=1$. Generalizing results in Section~\ref{s:cyclic}, we show
that the cyclic $q,t$-derangements satisfy natural $q,t$-analogs of
\eqref{e:derange}, \eqref{e:recur} and \eqref{e:recur1}. These results
also generalize formulas of Garsia and Remmel \cite{GarsiaRemmel1980}
who first introduced $q$-analogs for \eqref{e:recur} and
\eqref{e:recur1} for $r=1$ using a different (though equi-distributed)
permutation statistic.

Brenti \cite{Brenti1989} gave another $q$-analog for derangements of
$\Sn$ $q$-counted by weak excedances and conjectured many nice
properties for these numbers that were later proved by Canfield
(unpublished) and Zhang \cite{Zhang1995}. More recently, Chow
\cite{Chow2009} and Chen, Tang and Zhao \cite{CTZ2009} independently
extended these results to derangements of the hyperoctahedral
group. In Section~\ref{s:qcyclic}, we show that these results are
special cases of cyclic derangements of $\CS$ $q$-counted by a
generalization of weak excedances.

In Section~\ref{s:open}, we discuss possible directions for further
study.

\section*{Acknowledgments}

The author thanks Persi Diaconis, Ira Gessel, Gary Gordon, Vic Reiner
and Jeff Remmel for helpful discussions and suggestions for
generalizations.

\section{Cyclic derangements}\label{s:cyclic} 

The wreath product $\CS$ is the semi-direct product $(\Cr)^{\times n}
\rtimes \Sn$, where the symmetric group $\Sn$ acts on $n$ copies of
the cyclic group $\Cr$ by permuting the coordinates. Let $\zeta$ be a
generator for $\Cr$, e.g. take $\zeta$ to be a primitive $r$th root of
unity. We regard an element $\sigma \in \CS$ as a word $\sigma =
(\zeta^{e_1} s_1, \ldots, \zeta^{e_n} s_n)$ where $e_i \in
\{0,\ldots,r-1\}$ and $\{s_1,\ldots,s_n\} = \{1,\ldots,n\}$. Observe
that $|\CS| = r^n n!$, since there are $n!$ choices for the underlying
permutation $(s_1,\ldots,s_n)$ and $r^n$ choices for the signs
$(e_1,\ldots,e_n)$.

In this section we show that all of the usual formulas and proofs for
classical derangement numbers generalize to these wreath products.  We
begin with (\ref{e:derange}), giving the following formula for the
number of cyclic derangements. The two proofs below are essentially
the same, though the first is slightly more direct while the latter
will be useful for establishing $q$ and $q,t$ analogs.

\begin{theorem}
  The number of cyclic derangements in $\CS$ is given by
  \begin{equation}
    d_n^{(r)} = r^n n! \sum_{i=0}^{n} \frac{(-1)^i}{r^i i!}.
    \label{e:derange-cyclic}
  \end{equation}
  \label{t:derange-cyclic}
\end{theorem}

\begin{proof}[Inclusion-Exclusion Proof]
  Let $A_i$ be the set of $\sigma \in \CS$ such that $\sigma_i = +1
  \cdot i$. Then $|A_{j_1} \cap \cdots \cap A_{j_i}| = r^{n-i}
  (n-i)!$, since the positions $j_1,\ldots,j_i$ are determined and the
  remaining $n-i$ positions may be chosen arbitrarily. Therefore by
  the Inclusion-Exclusion formula, we have
  \begin{eqnarray*}
    \left| \Dnr \right| 
    & = & \left| \CS \right| - \left| A_1 \cup \cdots \cup A_n \right|
    \\
    & = & r^n n! -  \sum_{i=1}^n \sum_{j_1 < \cdots < j_i} (-1)^{i-1}
    \left| A_{j_1} \cap \cdots \cap A_{j_i} \right| \\
    & = & \sum_{i=0}^n \binom{n}{i} (-1)^i r^{n-i} (n-i)! = r^n n!
    \sum_{i=0}^{n} \frac{(-1)^i}{r^i i!}.
  \end{eqnarray*}
\end{proof}

\begin{proof}[M\"{o}bius Inversion Proof]
  For $S = \{ s_1 < s_2 < \cdots < s_m \} \subseteq [n]$ and $\sigma
  \in \Cr \wr \mathfrak{S}_A$, define the \textit{reduction} of
  $\sigma$ to be the permutation in $\Cr \wr \mathfrak{S}_m$ that
  replaces $\zeta^{e_i}s_i$ with $\zeta^{e_i}i$. If $\sigma \in \CS$
  has exactly $k$ fixed points, then define $\der(\sigma) \in
  \mathfrak{D}_{n-k}^{(r)}$ to be the reduction of $\sigma$ to the
  non-fixed points. For example, $\der(5314762)=$reduction of
  $53172=43152$ and any signs are carried over.

  The map $\der$ is easily seen to be a $\binom{n}{k}$ to $1$ mapping
  of subset of cyclic permutations with exactly $k$ fixed points onto
  $\mathfrak{D}_{n-k}^{(r)}$. Therefore
  \begin{equation}
    r^n n! = \sum_{k=0}^{n} \binom{n}{k} d_{n-k}^{(r)}.
    \label{e:perm_der}
  \end{equation}
  The theorem now follows by M\"{o}bius inversion \cite{EC1}.
\end{proof}

An immediate consequence of Theorem~\ref{t:derange-cyclic} is that the
probability that a random element of $\CS$ is a derangement is
approximately (and very nearly) $e^{-1/r}$. This verifies the
intuition that as $n$ and $r$ grow, most elements of $\CS$ are in fact
derangements. Table~\ref{t:dnr} gives values for $d_n^{(r)}$ for $r
\leq 5$ and $n \leq 6$.

\begin{table}[ht]
  \begin{center}
    \caption{\label{t:dnr} Cyclic derangement numbers $d_n^{(r)}$ for
      $r,n \leq 5$.} 
    \begin{tabular}{c|rrrrrrr}
      $r\setminus n$ & 0 & 1 & 2 & 3 & 4 & 5 & 6 \\ \hline
      1 & 1 & 0 &  1 &   2 &     9 &     44 & 265 \\
      2 & 1 & 1 &  5 &  29 &   233 &   2329 & 27949 \\
      3 & 1 & 2 & 12 & 116 &  1393 &  20894 & 376093 \\
      4 & 1 & 3 & 25 & 299 &  4785 &  95699 & 2296777 \\
      5 & 1 & 4 & 41 & 614 & 12281 & 307024 & 9210721 
    \end{tabular}
  \end{center}
\end{table}

We also have the following generalization of
\cite{GoMc2009}(Proposition 3.2), giving a formula relating the number
of cyclic derangements with the number of permutation derangements.

\begin{proposition}
  For $r \geq 2$ we have 
  \begin{equation}
    d_n^{(r)} = \sum_{i=0}^{n} \binom{n}{i} r^i (r-1)^{n-i} d_i 
    \label{e:both-derange}
  \end{equation}
  where $d_i = | \mathfrak{D}_i |$ is the number of derangements in
  $\mathfrak{S}_i$.
\label{p:both-derange}
\end{proposition}

\begin{proof}
  The number of derangements $\sigma \in \CS$ with precisely $i$
  indices $j$ for which $|\sigma_j| \neq j$ is equal to $d_i r^i$
  (choose a permutation derangement of these indices and a sign for
  each) times $(r-1)^{n-i}$ (choose a nonzero sign for indices $k$
  such that $|\sigma_k| = k$).
\end{proof}

Gordon and McMahon \cite{GoMc2009} observed that for $r=2$, the
expression in \eqref{e:both-derange} is precisely the \textit{rising
  $2$-binomial transform} of the permutation derangement numbers as
defined by Spivey and Steil \cite{SpSt2006}. In general, this formula
gives an interpretation for the mixed rising $r$-binomial transform
and falling $(r-1)$-binomial transform of the permutation derangements
numbers.

The following two term recurrence relation for cyclic derangements
generalizes \eqref{e:recur}. We give two proofs of this recurrence,
one generalizing the classical combinatorial proof of \eqref{e:recur}
and the other using the exponential generating function for cyclic
derangements.

\begin{theorem}
  For $n \geq 2$, the number of cyclic derangements satisfy
  \begin{equation}
    d_n^{(r)} = (rn-1) d_{n-1}^{(r)} + r(n-1) d_{n-2}^{(r)},
    \label{e:recurr-cyclic}
  \end{equation}
  with initial conditions $d_0^{(r)}=1$ and $d_1^{(r)}=r-1$.
  \label{t:recurr-cyclic}
\end{theorem}

\begin{proof}[Combinatorial Proof]
  For $\sigma \in \Dnr$, consider the cycle decomposition of
  underlying permutation $|\sigma| \in \Sn$. There are three cases to
  consider. Firstly, if $n$ is in a cycle of length one, then there
  are $r-1$ choices for $\sigma(n) = \zeta^{e} n$ with $e > 0$ and
  $d_{n-1}^{(r)}$ choices for a derangement of the remaining $n-1$
  letters. If $n$ is in a cycle of length two in $|\sigma|$, then
  there are $r(n-1)$ choices for the other occupant of the cycle in
  $\sigma$ and $d_{n-2}^{(r)}$ choices for a cyclic derangement of the
  remaining $n-2$ letters. Finally, if $n$ is in a cycle of length
  three or more, then there are $n-1$ possible positions for $n$ in
  $|\sigma|$, $r$ choices for the sign of $\sigma(n)$ and
  $d_{n-1}^{(r)}$ choices for a derangement of the remaining $n-1$
  letters. Combining these cases, we have
  \begin{displaymath}
    d_n^{(r)} = (r-1) d_{n-1}^{(r)} + r(n-1) d_{n-2}^{(r)} + r(n-1)
    d_{n-1}^{(r)}, 
  \end{displaymath}
  from which \eqref{e:recurr-cyclic} now follows.
\end{proof}

\begin{proof}[Algebraic Proof]
  First note that for fixed $r$,
  \begin{eqnarray*}
    \frac{e^{-x}}{1-rx} & = & \left( \sum_{i \geq 0} \frac{(-1)^i}{i!} x^i \right)
    \left(\sum_{j \geq 0} r^j x^j \right) \\
    & = & \sum_{n \geq 0} \sum_{i+j=n} \binom{n}{i} \frac{(-1)^i r^j}{i!}
    x^n = \sum_{n \geq 0} d_n^{(r)} \frac{x^n}{n!}
  \end{eqnarray*}
  is the exponential generating function for the number of cyclic
  derangements. Denoting this function by $D^{(r)}(x)$, we compute
  \[\begin{array}{l}
    \displaystyle \sum \left((rn-1) d^{(r)}_{n-1} + (rn-r)
      d^{(r)}_{n-2}\right) \frac{x^n}{n!} \\ [2ex]
    \displaystyle \hspace{1em} 
    = r \sum d^{(r)}_{n-1} \frac{x^n}{(n-1)!}
    - \sum d^{(r)}_{n-1} \frac{x^n}{n!}  
    + r \sum d^{(r)}_{n-2} \frac{x^n}{(n-1)!} 
    - r \sum d^{(r)}_{n-2} \frac{x^n}{n!} \\ [2ex]
    \displaystyle \hspace{1em} 
    = r x D^{(r)}(x)
    - \!\int\! D^{(r)}(x)
    + r x \!\int\! D^{(r)}(x)
    - r \!\int\!\!\int\! D^{(r)}(x) = D^{(r)}(x),
  \end{array}\]
  from which the recurrence now follows. 
\end{proof}

The following single term recurrence relation was discovered by Gordon
and McMahon, generalizing their result for the case $r=2$. 

\begin{corollary}
  For $n \geq 1$, the number of cyclic derangements satisfy
  \begin{equation}
    d_n^{(r)} = rn d_{n-1}^{(r)} + (-1)^n,
    \label{e:recur-cyclic}
  \end{equation}
  with initial condition $d_0^{(r)}=1$.
  \label{c:recur-cyclic}
\end{corollary}

This recurrence follows by induction from the formula in
Theorem~\ref{t:derange-cyclic} or the two term recurrence in
Theorem~\ref{t:recurr-cyclic}, though it would be nice to have a
direct combinatorial proof similar to that of Remmel \cite{Remmel1983}
for the case $r=1$.

\section{Cyclic $q,t$-derangements by major index}\label{s:qtcyclic} 

Gessel \cite{GeRe1993} derived a $q$-analog for the number of
permutation derangements as a corollary to a generating function
formula for counting permutations in $\Sn$ by descents, major index
and cycle structure. In order to state Gessel's formula, we begin by
recalling the $q$-analog of a positive integer $i$ given by $[i]_q = 1
+ q + \cdots + q^{i-1}$. In the same vein, we also have $[i]_q! =
[i]_q [i-1]_q \cdots [1]_q$, where $[0]_q!$ is defined to be $1$.

For a permutation $\sigma \in \Sn$, the \textit{descent set of
  $\sigma$}, denoted by $\Des(\sigma)$, is given by $\Des(\sigma) = \{
i \ | \ \sigma(i) > \sigma(i+1) \}$. MacMahon \cite{MacMahon1913} used
the descent set to define a fundamental permutation statistic, called
the \textit{major index} and denoted by $\maj(\sigma)$, given by
$\maj(\sigma) = \sum_{i \in \Des(\sigma)} i$.  Finally, recall
MacMahon's formula \cite{MacMahon1913} for $q$-counting permutations
by the major index statistic
\[ \sum_{\sigma \in \Sn} q^{\maj(\sigma)} = [n]_q! . \] 
Along these lines, define the \textit{$q$-derangement numbers},
denoted by $d_n(q)$, by
\begin{equation}
  d_n(q) = \sum_{\sigma \in \Dn} q^{\maj(\sigma)}.
\label{e:q-derange}
\end{equation}
Gessel showed that the $q$-derangement numbers for $\Sn$ are given by
\begin{equation}
  d_n(q) = [n]_q! \sum_{i=0}^{n} \frac{(-1)^i}{[i]_q!} q^{\binom{i}{2}}.
  \label{e:q-derange2}
\end{equation}

A nice bijective proof of (\ref{e:q-derange2}) is given by Wachs in
\cite{Wachs1989}, where she constructs a descent-preserving bijection
between permutations with specified derangement positions and shuffles
of two permutations and then makes use of a formula of Garsia and
Gessel \cite{GaGe1979} for $q$-counting shuffles. Garsia and Remmel
\cite{GarsiaRemmel1980} also studied $q$-derangement numbers using the
inversion statistic which is known to be equi-distributed with major
index.

Gessel's formula was generalized to the hyperoctahedral group by Chow
\cite{Chow2006} using the flag major index statistic. Here, we
generalize this further to $\CS$, and while the formula we derive
specializes to Chow's in the case $r=2$ and $t=q$, the major index
statistic we use differs. 

We begin with a generalized notion of descents derived from the
following total order on elements of $(\Cr \times [n]) \cup \{0\}$:
\begin{equation}
   \zeta^{r-1} n < \cdots < \zeta n < \zeta^{r-1} (n\!-\!1) < \cdots < 
   \zeta 1 < 0 < 1 < 2 < \cdots < n 
\label{e:total1}
\end{equation}
For $\sigma \in \CS$, an index $0 \leq i < n$ is a \textit{descent} of
$\sigma$ if $\sigma_i > \sigma_{i+1}$ with respect to this total
ordering, where we set $\sigma_0=0$. Note that for $\sigma \in \Sn$,
this definition agrees with the classical one. As with
permutations, define the \textit{major index} of $\sigma$ by
$\maj(\sigma) = \sum_{i \in \Des(\sigma)} i$. We also want to track
the signs of the letters of $\sigma$, which we do with the statistic
$\sign(\sigma)$ defined by $\sign(\sigma) = e_1 + \cdots + e_n$, where
$\sigma = (\zeta^{e_1} s_1, \ldots, \zeta^{e_n} s_n)$ . This is a
generalization of the same statistic introduced by Reiner in
\cite{Reiner1993}.

A first test that these statistics are indeed natural is to see that
the $q,t$ enumeration of elements of $\CS$ by the major index, sign
gives
\begin{equation}
  \sum_{\sigma \in \CS} q^{\maj(\sigma)} t^{\sign(\sigma)} = [r]_t^n [n]_q! ,
  \label{e:cyclic-qt}
\end{equation}
which is the natural $(q,t)$-analog for $r^n n! = |\CS|$. 

Analogous to (\ref{e:q-derange}), define the cyclic
$(q,t)$-derangement numbers by
\begin{equation}
  d_n^{(r)}(q,t) = \sum_{\sigma \in \Dnr}
  q^{\maj(\sigma)} t^{\sign(\sigma)}.
\label{e:cyclic-qt-derange}
\end{equation}
In particular, $d_n^{(1)}(q,t) = d_n(q)$ as defined in
\eqref{e:q-derange}. In general, we have the following $(q,t)$-analog
of \eqref{e:derange-cyclic} that specializes to \eqref{e:q-derange2}
when $r=1$ and to Chow's formula \cite{Chow2006}(Theorem 5)
when $r=2$ and $t=q$.

\begin{theorem}
  The cyclic $(q,t)$-derangement numbers are given by
  \begin{equation}
    d_n^{(r)}(q,t) = [r]_t^n [n]_q! \sum_{i=0}^{n}
    \frac{(-1)^i}{[r]_t^i [i]_q!} q^{\binom{i}{2}}. 
    \label{e:cyclic-q-derange}
  \end{equation}
  \label{t:cyclic-q-derange}
\end{theorem}

The proof of Theorem~\ref{t:cyclic-q-derange} is completely analogous
to Wachs's proof \cite{Wachs1989} for $\Sn$ which generalizes the
second proof of Theorem~\ref{t:derange-cyclic}. To begin, we define a
map $\varphi$ that is a sort of inverse to the map $\der$. Say that
$\sigma_i$ is a \textit{subcedant} of $\sigma$ if $\sigma_i < i$ with
respect to the total order in \eqref{e:total1}, and let $\sub(\sigma)$
denote the number of subcedants of $\sigma$.  For $\sigma \in \Cr \wr
\mathfrak{S}_m$, let $s_1 < \cdots < s_{\sub(\sigma)}$ be the absolute
values of subcedants of $\sigma$. If $\sigma$ has $k$ fixed points,
let $f_1 < \cdots < f_k$ be the (absolute values of) fixed points of
$\sigma$. Finally, let $x_1 > \cdots > x_{m-\sub(\sigma)-k}$ be the
absolute values of the remaining letters in $[m]$. For fixed $n$,
$\varphi(\sigma)$ is obtained from $\sigma$ by the following
replacements: 
\[\zeta^{e_i}s_i \mapsto \zeta^{e_i}i \hspace{2em} \zeta^{e_i}f_i
\mapsto \zeta^{e_i}(i+\sub(\sigma)) \hspace{2em} \zeta^{e_i}x_i
\mapsto \zeta^{e_i}(n-i+1).\]
For example, for $n=8$ and $\sigma = 326541$, we have $\varphi(\sigma)
= 638721$, and any signs are carried over.

Recall that for disjoint sets $A$ and $B$, a \textit{shuffle} of
$\alpha \in \Cr \wr \mathfrak{S}_A$ and $\beta \in \Cr \wr
\mathfrak{S}_B$ is an element of $\Cr \wr \mathfrak{S}_{A \cup B}$
containing $\alpha$ and $\beta$ as complementary subwords. Let
$\mathrm{Sh}(\alpha,\beta)$ denote the set of shuffles of $\alpha$ and
$\beta$. Then we have the following generalization of
\cite{Wachs1989}(Theorem 2).

\begin{lemma}
  Let $\alpha \in \mathfrak{D}_{n-k}^{(r)}$ and $\gamma =
  (\sub(\alpha)+1,,\ldots,\sub(\alpha)+k)$. Then the map $\varphi$
  gives a bijection $\{ \sigma \in \CS \ | \ \der(\sigma) = \alpha \}
  \stackrel{\sim}{\longrightarrow} \mathrm{Sh}(\varphi(\alpha),
  \gamma)$ such that $\Des(\varphi(\sigma)) = \Des(\sigma)$ and
  $\sign(\varphi(\sigma)) = \sign(\sigma)$.
  \label{l:preserve}
\end{lemma}

\begin{proof}
  The preservation of $\sign$ is obvious by construction. To see that
  the descent set is preserved, note that $\sign(\sigma_i)>0$ only if
  $\sigma_i$ is a subcedant and the relative order of subcedants,
  fixed points and the remaining letters is preserved by the map. It
  remains only to show that $\varphi$ is an invertible map with image
  $\mathrm{Sh}(\varphi(\alpha), \gamma)$. For this, the proof of
  \cite{Wachs1989}(Theorem 2) carries through verbatim thanks to the
  total ordering in \eqref{e:total1}.
\end{proof}

The only remaining ingredient to prove
Theorem~\ref{t:cyclic-q-derange} is the formula of Garsia and Gessel
\cite{GaGe1979} for $q$-counting shuffles. Though their theorem was
stated only for $\Sn$, using the total ordering in \eqref{e:total1}
the result holds in this more general setting.

\begin{lemma}
  Let $\alpha$ and $\beta$ be cyclic permutations of lengths $a$ and
  $b$, respectively, and let $\mathrm{Sh}(\alpha,\beta)$ denotes the
  set of shuffles of $\alpha$ and $\beta$. Then
  \begin{equation}
    \sum_{\sigma \in \mathrm{Sh}(\alpha,\beta)} q^{\maj(\sigma)} t^{\sign(\sigma)} 
    = {a+b \brack a}_q
    q^{\maj(\alpha) + \maj(\beta)} t^{\sign(\alpha) + \sign(\beta)} .
    \label{e:shuffle}
  \end{equation}
\label{l:shuffle}
\end{lemma}

\begin{proof}[Proof of Theorem~\ref{t:cyclic-q-derange}]
  For $\gamma$ as in Lemma~\ref{l:preserve}, observe $\maj(\gamma) = 0
  = \sign(\gamma)$. Thus applying Lemma~\ref{l:preserve} followed by
  Lemma~\ref{l:shuffle} allows us to compute
  \begin{eqnarray*}
    [r]_t^n [n]_q! & = & \sum_{\sigma \in \CS} q^{\maj(\sigma)}
    t^{\sign(\sigma)} \\
    & = & \sum_{k=0}^{n} \sum_{\alpha \in \mathfrak{D}_{n-k}^{(r)}}
    \sum_{\der(\sigma) = \alpha} q^{\maj(\sigma)} t^{\sign(\sigma)} \\
    & = & \sum_{k=0}^{n} \sum_{\alpha \in \mathfrak{D}_{n-k}^{(r)}}
    \sum_{\sigma \in \mathrm{Sh}(\varphi(\alpha),\gamma)}
    q^{\maj(\sigma)} t^{\sign(\sigma)} \\ 
    & = & \sum_{k=0}^{n} \sum_{\alpha \in \mathfrak{D}_{n-k}^{(r)}}
    {n \brack  k}_q q^{\maj(\alpha)} t^{\sign(\alpha)} \\
    & = & \sum_{k=0}^{n} {n \brack k}_q d_k^{(r)}(q,t).
  \end{eqnarray*}
  Applying M\"{o}bius inversion to the resulting equation yields
  \eqref{e:cyclic-q-derange}.
\end{proof}

Despite the ease with which Theorem~\ref{t:derange-cyclic} and the
proof carry through to this setting, there is no known expression for
$d_{n}^{(r)}(q,t)$ in terms of $d_{n}(q,t)$. Indeed, based on the
simple case with $r=2$ and $n=2$, no formula of the form of
Proposition~\ref{p:both-derange} exists.

The recurrence relation \eqref{e:recurr-cyclic} in
Theorem~\ref{t:recurr-cyclic} does have a natural $(q,t)$-analog. Note
that this specializes to the formula of Garsia and Remmel
\cite{GarsiaRemmel1980} in the case $r=1$. The proof we give is
combinatorial, though it would be nice to have a generating function
proof as well.

\begin{theorem}
  The cyclic $(q,t)$-derangement numbers satisfy 
  \begin{equation}
    d_n^{(r)}(q,t) = \left( [r]_t [n]_q - q^{n-1} \right)
    d_{n-1}^{(r)}(q,t) + \left( q^{n-1}[r]_t [n-1]_q \right)
    d_{n-2}^{(r)}(q,t), 
    \label{e:qtrecur}
  \end{equation}
  with initial conditions $d_0^{(r)}(q,t)=1$ and
  $d_1^{(r)}(q,t)=[r]_t-1$.
  \label{t:qtrecur}
\end{theorem}

\begin{proof}
  As in the combinatorial proof of Theorem~\ref{t:recurr-cyclic},
  consider the cycle decomposition of underlying permutation $|\sigma|
  \in \Sn$. We consider the same three cases, this time tracking the
  major index and sign. If $n$ is in a cycle of length one, then the
  $r-1$ choices for $\sigma(n) = \zeta^{e} n$ with $e > 0$ contribute
  $t[r-1]_t$, and there will necessarily be a descent in position
  $n-1$, thus contributing $q^{n-1}$. This case then contributes
  \[ t[r-1]_t q^{n-1} d_{n-1}^{(r)}(q,t). \]

  If $n$ is in a cycle of length two in $|\sigma|$, then the sign
  contribution of $\sigma_n$ is arbitrary contributing $[r]_t$. The
  $n-1$ choices for the other occupant of the cycle will add at least
  $n-1$ to the major index beyond the major index of the permutation
  with these two letters removed. This contributes a term of
  $q^{n-1}[n-1]_q$, making the total contribution
  \[ [r]_t q^{n-1} [n-1]_q d_{n-2}^{(r)}(q,t). \]

  Finally, if $n$ is in a cycle of length three or more, then each of
  the $n-1$ possible positions for $n$ in $|\sigma|$ increases the
  major index by one, contributing a factor of $[n-1]_q$. The $r$
  choices for the sign of $\sigma(n)$ again contribute $[r]_t$, and so
  we have a total of
  \[ [r]_t [n-1]_q d_{n-1}^{(r)}(q,t). \] 
  Adding these three cases yields \eqref{e:qtrecur}.
\end{proof}

As before, we may use induction and the above relation to derive the
following single term recurrence relation for cyclic
$(q,t)$-derangements generalizing \eqref{e:recur-cyclic} of
Corollary~\ref{c:recur-cyclic}.

\begin{corollary}
  The cyclic $(q,t)$-derangement numbers satisfy 
  \begin{equation}
    d_n^{(r)} = [r]_t [n]_q d_{n-1}^{(r)} + (-1)^nq^{\binom{n}{2}},
    \label{e:qtrecur1}
  \end{equation}
  with initial condition $d_0^{(r)}(q,t)=1$.
\end{corollary}

\begin{remark}
  There is another total ordering on elements of $(\Cr \times [n])
  \cup \{0\}$ that is equally as natural as the order given in
  \eqref{e:total1}, namely
  \begin{equation}
    \zeta^{r-1} n < \cdots < \zeta^{r-1} 1 < \zeta^{r-2} n < \cdots < 
    \zeta 1 < 0 < 1 < 2 < \cdots < n.
  \end{equation}
  While using this alternate order will result in a different descent
  set and major index for a given element of $\CS$, the distribution
  of descent sets over $\CS$ and even $\Dnr$ is the same with either
  ordering. In fact, there are many possible total orderings that
  refine the ordering on positive integers and yield the same
  distribution over $\CS$ and $\Dnr$, since the proof of
  Theorem~\ref{t:cyclic-q-derange} carries through easily for these
  orderings as well. We have chosen to work with the ordering in
  \eqref{e:total1} primarily to facilitate the combinatorial proof of  
  Theorem~\ref{t:qtrecur}.
\end{remark}

\section{Cyclic $q$-derangements by weak excedances}\label{s:qcyclic} 

Brenti \cite{Brenti1989} studied a different $q$-analog of derangement
numbers, defined by $q$-counting derangements by the number of
weak excedances, in order to study certain symmetric functions introduced
by Stanley \cite{Stanley1989}. Later, Brenti \cite{Brenti1994} defined
weak excedances for the signed permutations to study analogous
functions for the hyperoctahedral group. Further results were
discovered by Zhang \cite{Zhang1995} and Chow \cite{Chow2009} and
Chen, Tang and Zhao \cite{CTZ2009} for the symmetric group and
hyperoctahedral group, respectively. Below we extend these results to
the wreath product $\CS$. 

Recall that an index $i$ is a \textit{weak excedant} of $\sigma$ if
$\sigma(i)=i$ or $\sigma^2(i) > \sigma(i)$. It has long been known that
the number of descents and the number of weak excedances are
equi-distributed over $\Sn$ and that both give the \textit{Eulerian
  polynomials} $A_n(q)$:
\begin{equation}
  A_n(q) \stackrel{\mathrm{def}}{=} \sum_{\sigma \in \Sn}
  q^{\des(\sigma)} = \sum_{\sigma \in \Sn} q^{\exc(\sigma)} ,
\label{e:eulerian}
\end{equation}
where $\des(\sigma)$ is the number of descents of $\sigma$ and
$\exc(\sigma)$ is the number of weak excedances of $\sigma$.  We
generalize these statistics to $\CS$ by saying $1 \leq i \leq n$ is a
\textit{weak excedant} of $\sigma \in \CS$ if $\sigma(i)=i$ or if
$|\sigma(i)| \neq i$ and $\sigma^2(i) > \sigma(i)$ with respect to the
total order in \eqref{e:total1}.

As with the number of descents, this statistic agrees with the
classical number of weak excedances for permutations and Brenti's
statistic for signed permutations. Moreover, the equi-distribution of
the number of non-descents and the number of weak excedances holds in
$\CS$, and the same bijective proof using canonical cycle form
\cite{EC1} holds in this setting. Therefore define the \textit{cyclic
  Eulerian polynomial} $A^{(r)}_n(q)$ by
\begin{equation}
  A^{(r)}_n(q) = \sum_{\sigma \in \CS} q^{n-\des(\sigma)} =
  \sum_{\sigma \in \CS} q^{\exc(\sigma)}.  
\label{e:eulerianCS}
\end{equation}
Note that $A^{(r)}_n(q)$ is palindromic for $r \leq 2$,
i.e. $A^{(r)}_n(q) = q^{n} A^{(r)}_n(1/q)$. In particular,
\eqref{e:eulerianCS} specializes to \eqref{e:eulerian} when $r=1$ and
to Brenti's type B Eulerian polynomial when $r=2$. For $r \geq 3$, the
cyclic Eulerian polynomial is not palindromic.

Restricting to the set of cyclic derangements of $\CS$, the number of
descents and weak excedances are no longer equi-distributed, even for
$r=1$. Define the cyclic $q$-derangement polynomials $D^{(r)}_n(q)$ by
\begin{equation}
  D^{(r)}_n(q) = \sum_{\sigma \in \Dnr} q^{\exc(\sigma)}.
\end{equation}
We justify this definition with the following two-term recurrence
relation generalizing Theorem~\ref{t:recurr-cyclic}. Note that this
reduces to the result of Brenti \cite{Brenti1989} when $r=1$ and the
analog for the hyperoctahedral group \cite{Chow2009,CTZ2009} when
$r=2$.

\begin{theorem}
  For $n \geq 2$, the cyclic $q$-derangement polynomials satisfy
  \begin{equation}
    D_n^{(r)} = (n-1)rq \left(D_{n-1}^{(r)} + D_{n-2}^{(r)} \right) +
    (r-1) D_{n-1}^{(r)} + rq(1-q) \left(D_{n-1}^{(r)}\right)^{\prime}
    \label{e:recur-cyclic-q}
  \end{equation}
  with initial conditions $D_0^{(r)}(q)=1$ and $D_1^{(r)}(q)=r-1$.
  \label{t:recurr-cyclic-q}
\end{theorem}

\begin{proof}
  As in the combinatorial proof of Theorem~\ref{t:recurr-cyclic},
  consider the cycle decomposition of underlying permutation $|\sigma|
  \in \Sn$. We consider the same three cases, now tracking the number
  of weak excedances.  If $n$ is in a cycle of length one, then there
  are $r-1$ choices for $\sigma(n) = \zeta^{e} n$ with $e > 0$, and
  $n$ is not a weak excedant of $\sigma$. Removing this cycle leaves
  a cyclic derangement in $\mathfrak{D}_{n-1}^{(r)}$ with the same
  number of weak excedances, thus contributing
  \[ (r-1) D_{n-1}^{(r)}(q). \]

  If $n$ is in a cycle of length two in $|\sigma|$, then the sign
  contribution of $\sigma_n$ is arbitrary and there are $n-1$ choices
  for the other occupant of the cycle, say $k$. Moreover, exactly one
  of $k$ and $n$ will be a weak excedant, depending on the choice of
  sign, and so the contribution in this case is
  \[ r(n-1)q D_{n-2}^{(r)}(q). \]

  Finally, if $n$ is in a cycle of length three or more, then the sign
  of $\sigma_n$ is arbitrary is again arbitrary, but the affect on
  weak excedances for the $n-1$ possible placements of $n$ is more
  subtle. If $n$ is placed between $i$ and $j$ with $\sigma(j) =
  \sigma(\sigma(i)) > \sigma(i)$, then the number of (weak) excedances
  remains unchanged when inserting $n$. However, if $n$ is placed
  between $i$ and $j$ with $\sigma(j) = \sigma(\sigma(i)) <
  \sigma(i)$, then the insertion of $n$ creates a new (weak)
  excedant. Therefore we may count these cases by
  \begin{eqnarray*}
    r  \sum_{\tau \in \mathfrak{D}_{n-1}^{(r)}} \left( \exc(\tau)
      q^{\exc(\tau)} + (n-1-\exc(\tau)) q^{\exc(\tau)+1} \right) && \\
    = r(n-1)q D_{n-1}^{(r)}(q) + r(1-q)q \left( D_{n-1}^{(r)}
    \right)^{\prime}(q) . && 
  \end{eqnarray*}
  Adding these three cases yields \eqref{e:recur-cyclic-q}.
\end{proof}

Using \eqref{e:recur-cyclic-q}, we can also compute the exponential
generating function of the cyclic Eulerian polynomials and cyclic
$q$-derangement polynomials.

\begin{proposition}
  For $r \geq 1$, we have
  \begin{equation}
    \sum_{n \geq 0} A^{(r)}_n(q) \frac{x^n}{n!} = 
    \frac{\displaystyle (1-q)e^{x(1-q)}}{\displaystyle 1 - qe^{rx(1-q)}},
    \label{e:AqGen}
  \end{equation}
  and
  \begin{equation}
    \sum_{n \geq 0} D^{(r)}_n(q) \frac{x^n}{n!} = 
    \frac{\displaystyle (1-q)e^{x(r-1)}}{\displaystyle e^{qrx} - qe^{rx}}.
    \label{e:DqGen}
  \end{equation}
\end{proposition}

\begin{proof}
  Reversing the generating function proof of
  Theorem~\ref{t:recurr-cyclic}, it is straightforward to show that
  \eqref{e:DqGen} satisfies the recurrence relation in
  \eqref{e:recur-cyclic-q}. Enumerating elements of $\CS$ by the
  number of fixed points yields
  \begin{equation}
    A_n^{(r)}(q) = \sum_{k=0}^{n} \binom{n}{k} q^k D_{n-k}^{(r)}(q),
    \label{e:AD}
  \end{equation}
  from which \eqref{e:AqGen} follows.
\end{proof}

Recall that a sequence $a_0,a_1,\ldots,a_m$ of real numbers is
\textit{unimodal} if for some $j$ we have $a_0 \leq a_1 \leq \cdots
\leq a_j \geq a_{j+1} \geq \cdots \geq a_m$. A sequence is
\textit{log-concave} if $a_i^2 \geq a_{i-1} a_{i+1}$ for all $i$. It
is not difficult to show that a log-concave sequence of positive
numbers is unimodal. More generally, a sequence is a \textit{P\'{o}lya
  frequency sequence} if every minor of the (infinite) matrix
$(a_{j-i})$ is nonnegative, where we take $a_k=0$ for $k<0$ and
$k>m$. P\'{o}lya frequency sequence arise often in combinatorics, and
one of the fundamental results concerning them is the following.

\begin{theorem}
  The roots of a polynomial $a_0 + a_1 x + \cdots + a_m x^m$ are all
  real and nonpositive if and only if the sequence
  $a_0,a_1,\ldots,a_m$ is a P\'{o}lya frequency sequence.
\label{t:polya}
\end{theorem}

Using \eqref{e:recur-cyclic-q} and Theorem~\ref{t:polya}, we will show
that for fixed $n$ and $r$, the sequence $a_m = \# \{ \sigma \in
\Dnr \ | \ \exc(\sigma) = m\}$ is a  P\'{o}lya frequency sequence. 

\begin{theorem}
  For $n \geq 2$, the roots of the cyclic $q$-derangement polynomial
  $D_{n}^{(r)}(q)$ interlace the roots of $D_{n+1}^{(r)}(q)$. In
  particular, they are distinct, nonpositive real numbers
\end{theorem}

\begin{proof}
  We proceed by induction on $n$. For $r=2,3$, we have
  \[ D_2^{(r)}(q) = r^2q + (r-1)^2 \hspace{1em} \mbox{and}
  \hspace{1em} D_3^{(r)}(q) = r^3q^2 + (4r-3)r^2q + (r-1)^3. \] Thus
  the root of $D_2^{(r)}(q)$ is $-(r-1)^2/r^2$, which is indeed real
  and nonpositive and lies between the two distinct negative real
  roots of $D_3^{(r)}(q)$. This demonstrates the base case, so assume
  the result for $n-1 \geq 2$.

  Let $q_1 < q_2 < \cdots < q_{n-2} < 0$ be the simple roots of
  $D_{n-1}^{(r)}(q)$. It is straightforward to show that
  $(D_{n-1}^{(r)})^{\prime}(q_i)$ has sign $(-1)^{n-i}$, and by
  induction, the sign of $D_{n-2}^{(r)}(q_i)$ is also $(-1)^{n-i}$ as
  the roots are interlaced. By \eqref{e:recur-cyclic-q}, we have
  \[ D_{n}^{(r)}(q_i) = (n-1)rq_i D_{n-2}^{(r)}(q_i) + rq_i(1-q_i)
  \left(D_{n-1}^{(r)}\right)^{\prime}(q_i), \] 

  from which it follows that the sign of $D_{n}^{(r)}(q_i)$ is
  $(-1)^{n-i-1}$. Noting as well that the leading term of
  $D_{n}^{(r)}(q)$ is $r^n$ and the constant term is $(r-1)^n$, it
  follows from the Intermediate Value Theorem that the are roots of
  $D_{n}^{(r)}(q)$ are interlaced with $-\infty,q_1,\ldots,q_{n-2},0$.
\end{proof}

\section{Further directions}\label{s:open} 

\subsection*{Binomial transforms}

Spivey and Steil \cite{SpSt2006} define two variants of the binomial
transform of a sequence: the \textit{rising $k$-binomial transform}
and the \textit{falling $k$-binomial transform} given by
\begin{equation}
  r_n = \sum_{i=0}^{n} \binom{n}{i} k^i a_i
  \hspace{1em} \mbox{and} \hspace{1em} 
  f_n = \sum_{i=0}^{n} \binom{n}{i} k^{n-i} a_i,
\end{equation}
respectively. Gordon and McMahon noted that the number of derangements
in the hyperoctahedral group gives the rising $2$-binomial transform
of the derangement numbers for $\Sn$. More generally,
Proposition~\ref{p:both-derange} shows that the cyclic derangement
numbers $d_n^{(r)}$ give a mixed version of the rising $r$-binomial
transform and falling $(r-1)$-binomial transform of $d_n$. This new
hybrid $k$-binomial transform may share many of the nice properties of
Spivey and Steil's transforms, including Hankel invariance and/or a
simple description of the change in the exponential generating
function. Further, it could be interesting to evaluate the expression
in (\ref{e:both-derange}) for negative or even non-integer values of
$k$. For instance, taking $k=1/2$ gives the \textit{binomial mean
  transform} which is of some interest.

\subsection*{Limiting distributions}

The explicit expression in (\ref{e:derange-cyclic}) immediately gives
approximations and asymptotics for the probability that a random
element of $\CS$ is a derangement. Moreover, these formulae can be
used to calculate the number of elements with a given number of fixed
points, and so, too, they can be used to calculate the expected number
of fixed points of a random element. For instance, in $\Sn$ it is
know that the number of fixed points of a random permutation has a
limiting Poisson distribution, and recent work by Diaconis, Fulman and
Guralnick \cite{DFG2008} has extended this to primitive actions of
$\Sn$. A natural extension would be use the combinatorics presented
here to consider the imprimitive action of $\Sn$ in $\CS$.

\subsection*{Orthogonal idempotents}

Another direction would be to generalize the work of Schocker
\cite{Schocker2003} where he uses derangement numbers to construct $n$
mutually orthogonal idempotents in Solomon's descent algebra for
$\Sn$. In doing so, he also discovers a new proof of Gessel's formula
for $q$-derangements. Thus extending these techniques to analogs of
the descent algebra for $\CS$ could also lead to new proofs of the
formulae in Section~\ref{s:qtcyclic}.

\subsection*{Symmetric unimodal polynomials}

In \cite{Brenti1990}, Brenti used symmetric functions to define
several new classes of symmetric unimodal polynomials. Brenti showed
there is an explicit connection between these polynomials and the
$q$-Eulerian polynomials and $q$-derangement polynomials of $\Sn$
counted by weak excedances. Brenti generalized much of this work to
the hyperoctahedral group, with additional results extended by Chow
\cite{Chow2009}. A natural question is to see if there exist analogs
of these results for $\CS$ involving the polynomials $A^{(r)}_n(q)$
and $D^{(r)}_n(q)$ studied in Section~\ref{s:qcyclic}.

\bibliography{derangements}

\begin{thebibliography}{10}

\bibitem{Brenti1989}
F.~Brenti.
\newblock Unimodal, log-concave and {P}\'olya frequency sequences in
  combinatorics.
\newblock {\em Mem. Amer. Math. Soc.}, 81(413):viii+106, 1989.

\bibitem{Brenti1990}
F.~Brenti.
\newblock Unimodal polynomials arising from symmetric functions.
\newblock {\em Proc. Amer. Math. Soc.}, 108(4):1133--1141, 1990.

\bibitem{Brenti1994}
F.~Brenti.
\newblock {$q$}-{E}ulerian polynomials arising from {C}oxeter groups.
\newblock {\em European J. Combin.}, 15(5):417--441, 1994.

\bibitem{CTZ2009}
W.~Y.~C. Chen, R.~L. Tang, and A.~F.~Y. Zhao.
\newblock Derangement polynomials and excedances of type {$B$}.
\newblock {\em Electron. J. Combin.}, 16(2, Special volume in honor of Anders
  Bjorner):Research Paper 15, 16, 2009.

\bibitem{Chow2006}
C.-O. Chow.
\newblock On derangement polynomials of type {$B$}.
\newblock {\em S\'em. Lothar. Combin.}, 55:Art. B55b, 6 pp. (electronic), 2006.

\bibitem{Chow2009}
C.-O. Chow.
\newblock On derangement polynomials of type {B}. {II}.
\newblock {\em J. Combin. Theory Ser. A}, 116(4):816--830, 2009.

\bibitem{DFG2008}
P.~Diaconis, J.~Fulman, and R.~Guralnick.
\newblock On fixed points of permutations.
\newblock {\em J. Algebraic Combin.}, 28(1):189--218, 2008.

\bibitem{GaGe1979}
A.~M. Garsia and I.~Gessel.
\newblock Permutation statistics and partitions.
\newblock {\em Adv. in Math.}, 31(3):288--305, 1979.

\bibitem{GarsiaRemmel1980}
A.~M. Garsia and J.~Remmel.
\newblock A combinatorial interpretation of {$q$}-derangement and
  {$q$}-{L}aguerre numbers.
\newblock {\em European J. Combin.}, 1(1):47--59, 1980.

\bibitem{GeRe1993}
I.~M. Gessel and C.~Reutenauer.
\newblock Counting permutations with given cycle structure and descent set.
\newblock {\em J. Combin. Theory Ser. A}, 64(2):189--215, 1993.

\bibitem{GoMc2009}
G.~Gordon and E.~McMahon.
\newblock Moving faces to other places: Facet derangements.
\newblock {\em Amer. Math. Monthly}.
\newblock To appear.

\bibitem{MacMahon1913}
P.~A. MacMahon.
\newblock The {I}ndices of {P}ermutations and the {D}erivation {T}herefrom of
  {F}unctions of a {S}ingle {V}ariable {A}ssociated with the {P}ermutations of
  any {A}ssemblage of {O}bjects.
\newblock {\em Amer. J. Math.}, 35(3):281--322, 1913.

\bibitem{Reiner1993}
V.~Reiner.
\newblock Signed permutation statistics.
\newblock {\em European J. Combin.}, 14(6):553--567, 1993.

\bibitem{Remmel1983}
J.~B. Remmel.
\newblock A note on a recursion for the number of derangements.
\newblock {\em European J. Combin.}, 4(4):371--374, 1983.

\bibitem{Schocker2003}
M.~Schocker.
\newblock Idempotents for derangement numbers.
\newblock {\em Discrete Math.}, 269(1-3):239--248, 2003.

\bibitem{SpSt2006}
M.~Z. Spivey and L.~L. Steil.
\newblock The {$k$}-binomial transforms and the {H}ankel transform.
\newblock {\em J. Integer Seq.}, 9(1):Article 06.1.1, 19 pp. (electronic),
  2006.

\bibitem{Stanley1989}
R.~P. Stanley.
\newblock Log-concave and unimodal sequences in algebra, combinatorics, and
  geometry.
\newblock In {\em Graph theory and its applications: {E}ast and {W}est
  ({J}inan, 1986)}, volume 576 of {\em Ann. New York Acad. Sci.}, pages
  500--535. New York Acad. Sci., New York, 1989.

\bibitem{EC1}
R.~P. Stanley.
\newblock {\em Enumerative combinatorics. {V}ol. 1}, volume~49 of {\em
  Cambridge Studies in Advanced Mathematics}.
\newblock Cambridge University Press, Cambridge, 1997.
\newblock With a foreword by Gian-Carlo Rota, Corrected reprint of the 1986
  original.

\bibitem{Wachs1989}
M.~L. Wachs.
\newblock On {$q$}-derangement numbers.
\newblock {\em Proc. Amer. Math. Soc.}, 106(1):273--278, 1989.

\bibitem{Zhang1995}
X.~Zhang.
\newblock On {$q$}-derangement polynomials.
\newblock In {\em Combinatorics and graph theory '95, {V}ol.\ 1 ({H}efei)},
  pages 462--465. World Sci. Publ., River Edge, NJ, 1995.

\end{thebibliography}
\bibliographystyle{abbrv}

\end{document}